# Carrick mat and further development of bipartite knots

Alina Pavlikova *

## Abstract

In this work we continue to develop the theory of bipartite knots, that starts from question raised in 1987 by J. Przytycki[2]. Problem 1.60 in Kirby's list of open problems in topology [1]

# Table of contents

# Introduction

## History of a question

The problem raised by Jozef Przytycki in 1987 is to investigate which knots possess a matched diagram. This question appears in the well-known collection "Open problems in topology" maintained by Rob Kirby [1], as part of Problem 1.60. More exactly, Conjecture 1(a) therein (belonging to Przytycki [2]) reads: "There are oriented knots without a matched diagram". This conjecture stayed open for 24 years, notwithstanding the effort of several excellent mathematicians, including its author and J. H. Conway [3]. My teacher S. Duzhin give a positive solution to the conjecture, that is, demonstrate that some knots, e.g. pretzel knot with parameters $(3, 3, -3)$, are not bipartite. Also he proved next statment:

**Theorem 1.** *Let a knot $K$ admit a matched diagram and its elementary ideal $E_k(K)$ for some $k$ does not coincide with the ring $\mathbb{Z}[t, t^{-1}]$. Then $E_k(K)$ does not contain the element $t + 1$.*

This theorem provides only the necessary conditions for the existence of a matched knot diagram. In particular, the ideal $E_2(8_{18})$ is nontrivial, but does not contain the element $t + 1$. Due to the fact that it was not possible to find its paired diagram for a long time, and also because the ideal $E_2(K)$ is trivial for all known paired knots, the question arose as to whether paired knots with a non-trivial second ideal $E_2(K)$. The question about matched diagram of $8_{18}$ was open.

\*\*UNIGE, Villa Battelle, 1227 Carouge, Suisse kussike(at)gmail.com

The aim of this work is to further study the properties of the ideal $E_2(K)$ for paired nodes. Its main results are

- Presenting a matched diagram of a knot $8_{18}$;

- Representation of a family of paired nodes with nontrivial ideal $E_2(K)$;

- Introduction of two operations on bipartite knots whose diagram admits a so-called support chord (see the definition further), which do not change the ideal $E_2(K)$ of such knots;

- New proof of the triviality of the ideal $E_2(K)$ for rational knots (see the definition in Section 4). Presenting a paired knot with a trivial second elementary ideal that is not rational.

# Basic definitions and constructions

## Knot, knot diagram, knot group

**Definition.** A knot is a smooth embedding of $S^1$ into $\mathbb{R}^3$. Two knots $K_1$ and $K_2$ are said to be equivalent if there exist orientation-preserving homeomorphisms $\psi : \mathbb{R}^3 \to \mathbb{R}^3$ and $\varphi : S^1 \to S^1$ such that the diagram is commutative.

$$\begin{array}{ccc} S^1 & \xrightarrow{\phi} & S^1 \\ \downarrow{\scriptstyle K_1} & & \downarrow{\scriptstyle K_2} \\ \mathbb{R}^3 & \xrightarrow{\psi} & \mathbb{R}^3 \end{array}$$

Fig. 1: knot equivalence

*Remark.* In what follows, the node $K : S^1 \to \mathbb{R}^3$ and its image $K(S^1)$ will be denoted by the same letter $K$, if this does not cause ambiguity in interpretation.

**Definition.** The group of a knot $K \subset \mathbb{R}^3$ is the fundamental group of the complement to its image $\pi_1(\mathbb{R}^3 \setminus K)$.

*Remark.* Let us denote knot complement $\mathbb{R}^3 \setminus K$ by $X$.

Let us choose some hyperplane $P$ in $\mathbb{R}^3$. Consider the projection $d : \mathbb{R}^3 \to P$. Let $K$ be a knot, then $\tilde{K} = d(K)$ is some plane curve, which is a projection of the knot $K$. The orientation of $\tilde{K}$ is naturally induced. Note that the curve $d(K)$ has self-intersection points. By choosing an equivalent $K$ knot, one can achieve that the projection has a finite number of self-intersection points and any point of the projection has at most two inverse images. At each double point of the projection $d(K)$, using a small discontinuity, we can mark which branch of the knot with respect to the map $d$ in the preimage passed above and which was below.

**Definition.** Such projection $d(K)$ of a knot $K$ called a regular plane diagram of the knot.

*Remark.* In what follows, a knot diagram will mean a regular plane diagram. A knot diagram uniquely define its equivalence class.

## Alexander module

As is known from cover theory, to every subgroup $H$ in $\pi_1(X)$ there is corresponding covering of the space $X$, the fundamental group of which is exactly $H$ Take the covering $\tilde{X}$ corresponding to the subgroup $[\pi_1(X), \pi_1(X)]$. The projection $p : \tilde{X} \to X$ induces an injective homomorphism $p_* : \pi_1(\tilde{X}) \to \pi_1(X)$ and $p_*\pi_1(\tilde{X}) = [\pi_1(X), \pi_1(X)]$ ◂ $\pi_1(\tilde{X})$. Since this covering is normal, the automorphism group of the covering $\tilde{X}$ is isomorphic to the factor $\frac{\pi_1(X)}{[\pi_1(X),\pi_1(X)]}$. By the Hurewicz theorem, this factor is isomorphic to $H_1(X) = \mathbf{Z}$ Let $t$ be the generator of $H_1(X)$. Consider the automorphism $\tilde{X}$ corresponding to $t$ and denote it by $\tau$.

The Abelian group $H_1(\tilde{X})$ is naturally endowed with the structure of a module over the group ring $H_1(X)$ $\cdot$ $\mathbf{Z}$, which is the ring of integer Laurent polynomials $\mathbf{Z}[t, t^{-1}]$ : let $f(t) = a_{-m}t^{-m} + \cdots + a_0 + \cdots + a_m t^m$. Then for $\alpha \in H_1(\tilde{X})$ we put $f(t)\alpha = a_{-m}\tau_*^{-m} + \cdots + a_0 + \cdots + a_m \tau_*^m$. For any knot $K$ this module which we will call $M_K$ is finitely generated

**Definition.** The $M_K$ module is called the Alexander module of $K$.

## Trivial ideals of a knot

Corepresentation of a group $\pi_1(X)$ let us to obtain the matrix closely related to Alexander matrix of a module $M_K$. This is done as follows.

Consider oriented plane diagram $D$ of a knot $K$. In the neighborhood every double point diagram looks as depicted on Pic.

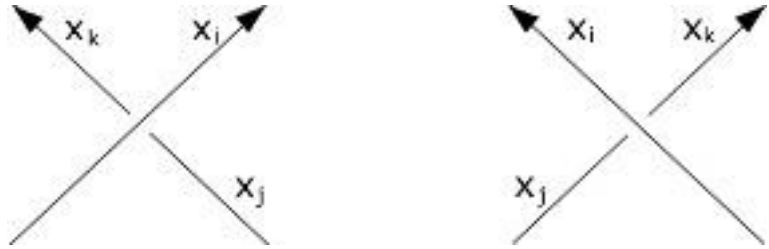

Fig. 2: Type of crossings. Letters denote arcs of diagram.

**Definition.** We call the connected component of $D$ an arc.

*Remark.* Note that number of crossings of diagram $D$ equal to the numbers of its arcs.

**Theorem 2.** *Every two Alexander matrix $A$ and $A'$ of module $M$ are equivalent.*

**Definition.** Let $A$ — Alexander matrix of module $M$ that have size $m \times n$. For every $k \in \mathbf{Z}$ define some ideal $E_k(M)$ of a ring $\mathbf{Z}[t, t^{-1}]$ by the next way:

- $0 < n - k \le m$, then $E_k(M)$ – ideal, generated by all minors of $A$ of order $n - k + 1$.

- $n - k > m$, then $E_k(M) = 0$.

- $n - k \le 0$, then $E_k(M) = \mathbf{Z}[t, t^{-1}]$ .

This ideal doesn't depend from choice of matrix $A$ in its equivalent class and called $k$-th elementary ideal (or $k$-th Alexander ideal) of module $M$.

**Definition.** For every $i$, $1 \le i \le m$ define the operator of free differential $\frac{\partial}{\partial \gamma_i} : F_m \to \mathbf{Z}[F_m]$,defined by the following identities:

1. $\frac{\partial ab}{\partial \gamma_i} = \frac{\partial a}{\partial \gamma_i} + a\frac{\partial b}{\partial \gamma_i}$

2. $\frac{\partial \gamma_i}{\partial \gamma_j} = \delta_{ij}$

Construct the matrix $(\frac{\partial r_i}{\partial x_j})$ and replace all the powers of the generators in it with the powers of the variable $t$: $x_i^k \to t^k$. The resulting matrix will be denoted by $Q_D$.

*Remark.* The $Q_D$ matrix is the Alexander matrix of the module $M_K \oplus \mathbf{Z}[t, t^{-1}]$.

The ideal generated by all minors of order $n - m$ of the matrix $Q_D$ coincides with the $m$-th elementary ideal of the module $M_K$.

## Alexander matrix and trivial ideals of the module over Laurent polynomial ring

**Definition.** Any finitely generated module over the ring $\mathbb{Z}[t, t^{-1}]$ can be represented as a quotient module of some free module. The $M_K$ module is called the Alexander module of $K$.

**Definition.** Alexander matrix $A$ of module $M$ over ring $\mathbb{Z}[t, t^{-1}]$ called matrix of the map $\beta$ between two modules $E$ and $F$ with basis' s $e_1 \ldots e_m$ and $f_1 \ldots f_n$, such that sequence is exact $E \xrightarrow{\beta} F \xrightarrow{\gamma} M \longrightarrow 0$

Every row $m$ from row of the matrix $A$ corresponds to the generator of $M$, and every column $n$ from columns – to the relation between generators.

## Bipartite knot and Chord diagram

**Definition.** Matched diagram of a knot is a diagram whose crossings are split in pairs of the types depicted in Fig.2.

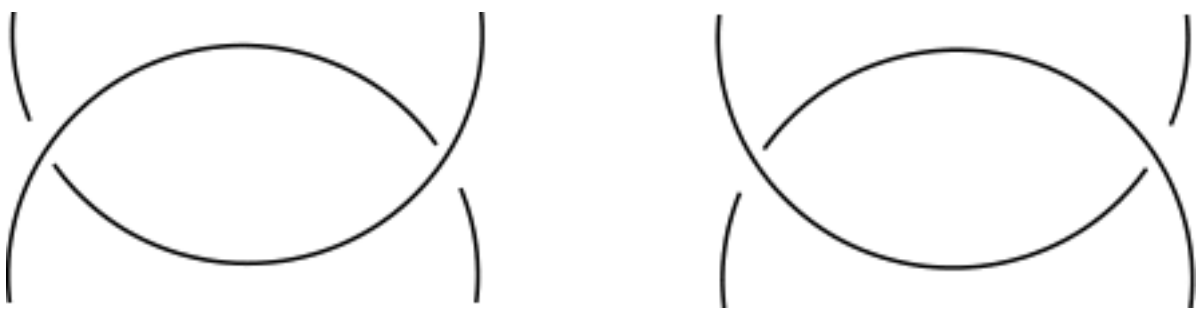

The left pair of crossings are said to be positive, those the right one, negative. We will call bipartite a knot that can be represented by a matched diagram.

**Definition.** Chord diagrm of order $n$ for $n \in \mathbb{Z}$ — its an oriented circle $S^1$ with matched $2n$ distinct points, divided by pairs named chords.

*Remark.* It' s easy to represent every chord by line segment, that connects corresponding points. Here $S^1$ considered under homeomorphisms preserved diffeomorphism orientation.

Every mathed diagram $D$ correspods chord diagram $C_D$. Namely, let $D$ — matched diagram. Change every pair of crossings to the pair of parallel segments, directed as corresponding part of knot, and connect it by common perpendicular.

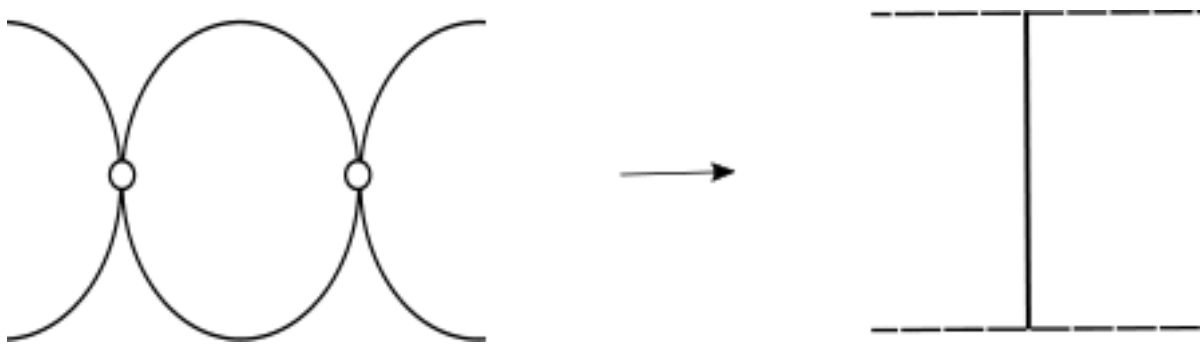

Fig. 3: Local transformation pair of crossings to a chord

**Example 1.** *On the Picture 4 we demonstrate chord diagram for matched diagram of trefoil.*

Arc of chord diagram called arc of the circle, they bijectively mapped to the arc of knot diagram. Here and throughout we think outer circle oriented counter-clockwise direction.

This procedure is reversible: from a bipartite signed chord diagram one can reconstruct the knot diagram in a unique way (see Fig. 1.5).

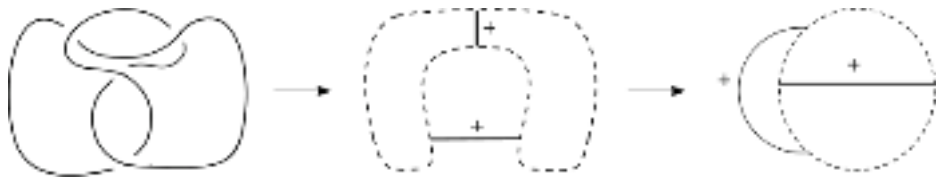

Fig. 4: trefoil

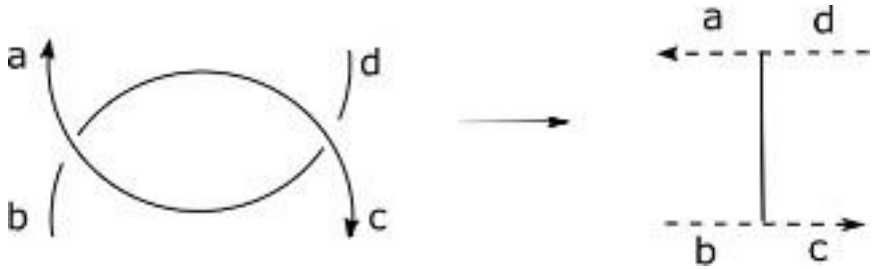

Fig. 5: arc corresponding

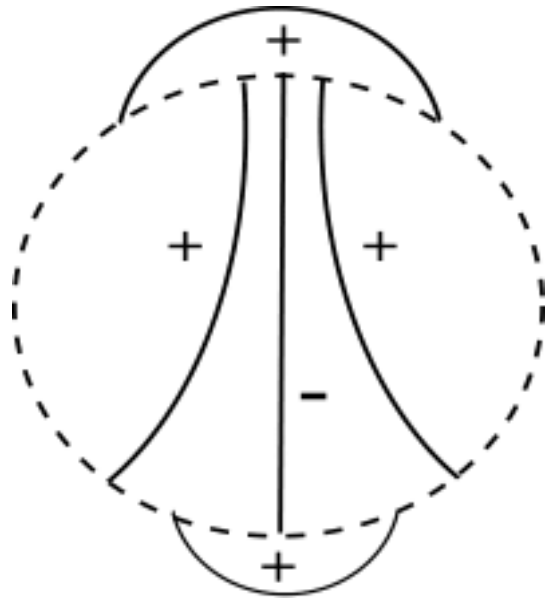

Fig. 6:

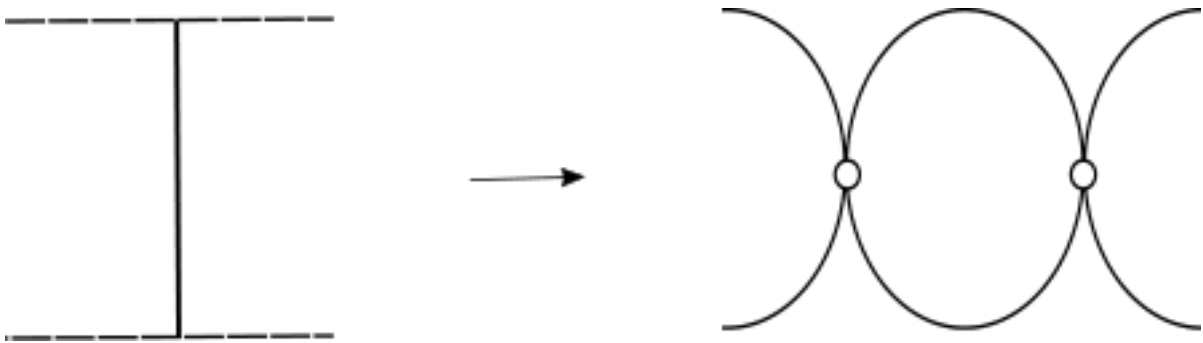

Fig. 7: chord to a pair of crossings

## $8_{18}$ **is bipartite!**

All knots up to 8 crossings, except for $8_{18}$, have a trivial second ideal and their paired diagrams were constructed by S.V. Duzhin. As the first example of a bipartite knot with a non-trivial second ideal, we give the knot $8_{18}$, its second elementary ideal $[t^2\_t + 1]$, the existence of its matched diagram was still an open question. It is shown in figure 8.

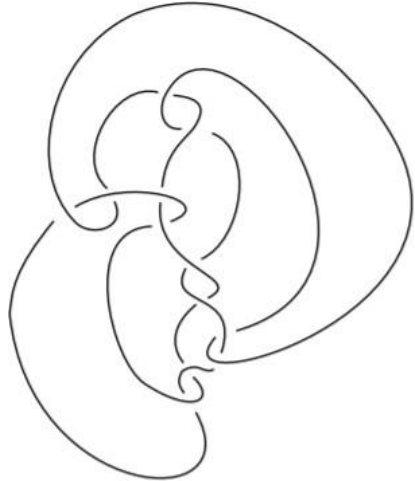

Fig. 8: matched diagram $8_{18}$.

## Knot union

The authors of the paper define an operation of the union of knots that generalise the operation of connected sum. Also in this article it is proved that union of two nontrivial knots results in a nontrivial knot.

*Remark.* In this part of the work we gonna use definitions and results from next papers We will consider only special case of union, we will take union of knot and its mirror.

Let $D$ be the diagram of some knot $K$. Consider a rectangle on the plane of the diagram of the knot $D$ Any finitely generated module over the ring $\mathbf{Z}[t, t^{-1}]$ can be represented as a quotient module of some free module. The $M_K$ module is called the Alexander module of $K$. end definition row and column turns out

Also in this article it is proved that pairing two nontrivial knots results in a nontrivial knot. special case

Let $D$ be the diagram of some knot $K$. Consider a rectangle on the plane of the diagram of the node $D$

such that all the crossings of the knot are inside this rectangle. Let us number the points of intersection of the knot with this rectangle according to their order when traversing the border of the rectangle clockwise (with some choice of the beginning of the traversal). We assume that the intersection of the knot diagram with the outside of the rectangle consists of a certain number of arcs $b_j$ $j \in (0, 1..., k)$ connecting pairs of neighboring intersection points of the knot diagram (according to the chosen numbering) with the boundary rectangle. Without loss of generality, we can assume that all points of intersection of the knot diagram with the boundary of the rectangle lie on one of its sides. Now consider the diagram $D^*$, obtained by reflecting $D$ and a rectangle about an axis parallel to the side of the rectangle that has a non-empty intersection with $D$. We choose the axis so that it is in the half-plane separated by a straight line that is an extension of the side of the rectangle in which all the arcs $b_j$ lie. The images of the arcs $b_j$ on the mirror-reflected diagram will be denoted, respectively, by $b_{\bar{j}}$.

Let' s take a set of integers $n_1, \ldots, n_j$. Replace all pairs of arcs, $(b_j, b_{\bar{j}})$ $j \in (1..., k)$, with special fragments with $n_j$ intersections, respectively, shown in the figure and unite the pair $(b_0, b_{\bar{0}})$ as shown in the figure

The result of this union with the parameters $n_1, \ldots, n_k$ will be called the duplicate of the note $K$ and denoted by $K_D(n_1, \ldots, n_k)$.

In the paper it is proved that the double $K_D(n_1, \ldots, n_k)$ of a nontrivial knot $K$ is nontrivial. We will be interested in the case when all numbers $n_1, \ldots, n_k$ are even.

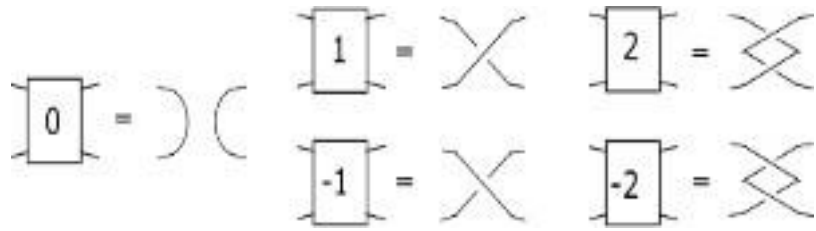

Fig. 9: Special fragments.

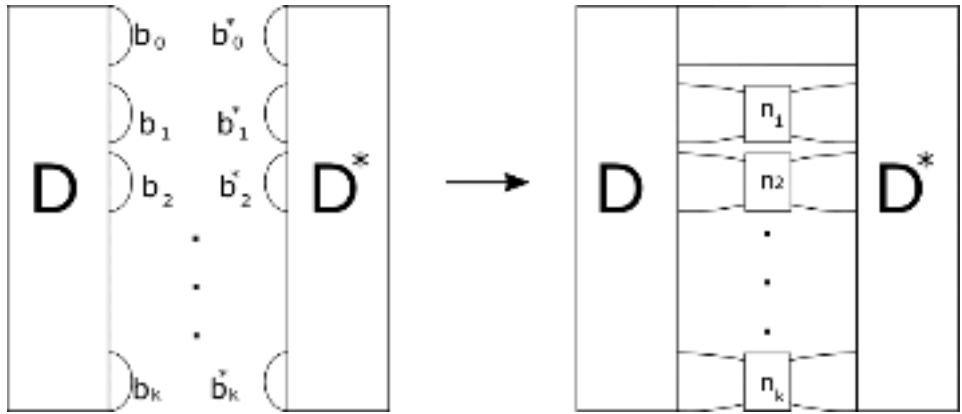

Fig. 10: Union of the diagram $D$ and $D*$.

**Example 2.** *Picture 11 illustrate one of the possible dublicate of trefoil.*

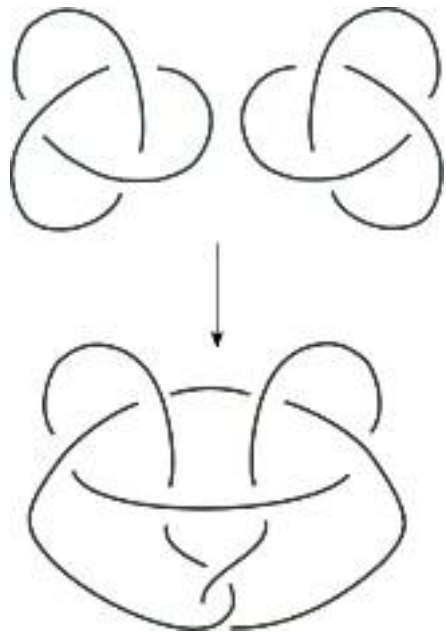

Fig. 11: trefoil duplicate with parameter $n_1 = 2$.

The family of duplicates of a trefoil generated by its standard diagram with $k = 1$ is called the Milnor-Fox family. The family of doubles of figure-eight knot generated by its standard diagram is called the Kanenobu family.

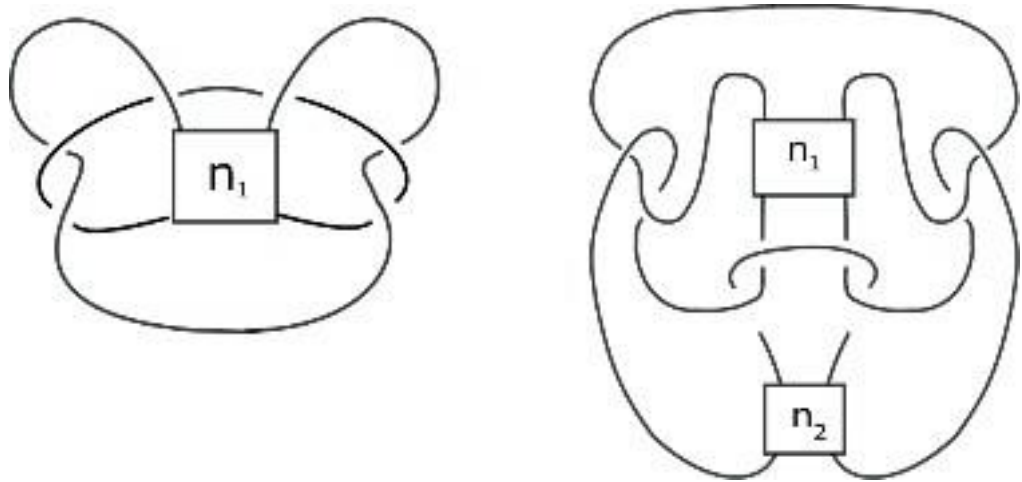

Fig. 12: Milnor-Fox and Kanenobu families

*Remark.* It's clear that the duplicate of the bipartite knot $K$ with its mirror image for even $n_j$ when choosing its mathced diagram will also be bipartite

For example, with Redemeister moves, knots from the Milnor-Fox family can be transformed in such way that their diagrams are matched

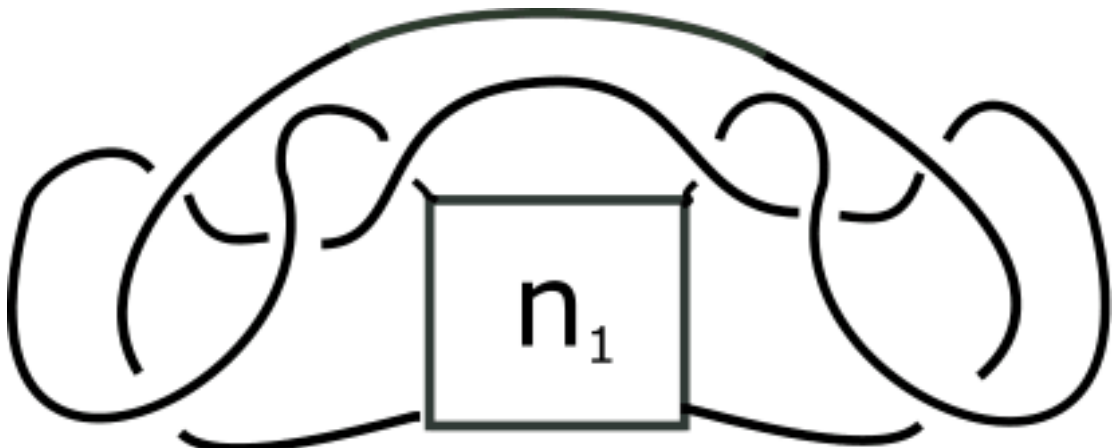

Fig. 13: Obtaining Milnor-Fox family with trefoil matched diagram

Family of chord diagram of such representation shown on the Fig

**Lemma 3.** *There exist infinitely many bipartite knots with nontrivial second Alexander ideal.*

*Proof.* It was shown in that knots that are duplicates with even coefficients have a special structure of the Alexander module. For example, for nodes from the Milnor-Fox family (See Fig.1.10 on the left) corresponding to an even $n_1$ the Alexander matrix, according to looks like :

$$\begin{matrix} t^2 - t + 1 & kt \\ 0 & t^2 - t + 1 \end{matrix}$$

where $2k = n_1$ is the number of intersections. The second trivial ideal is generated by the elements $[t^2 - t + 1, k]$. And this means that there are infinitely many bipartite knots with a nontrivial second elementary ideal. For example, in the case $k = 2$, we get a bipartite knot $10_{140}$ with a nontrivial second ideal. $\qquad\square$

**Lemma 4.** *There exist infinitely many bipartite knots with trivial second Alexander ideal.*

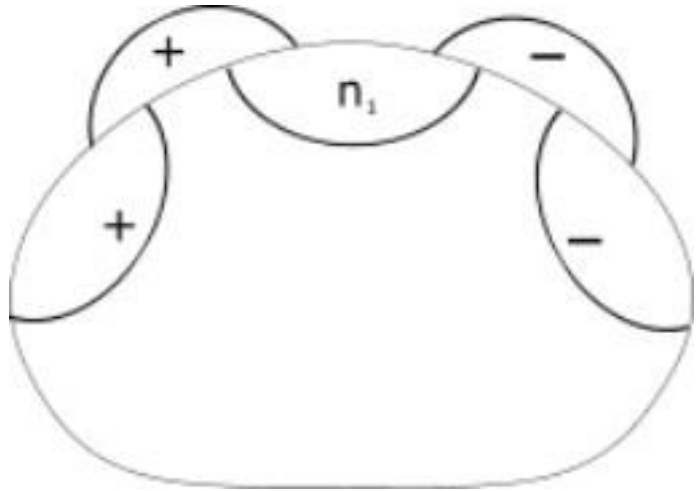

Fig. 14: Family of chord diagram for knots from Milnor-Fox family. Middle inner chord with index $n_1$ denote set of $|n_1|$ parallel chords with the same sign as $n_1$.

*Proof.* Note that one can make a duplicate of any knot from the Milnor-Fox family. The result is a new family of knots $\textbf{Sem}_1$. Family of chord diagram of such representation shown on the Fig

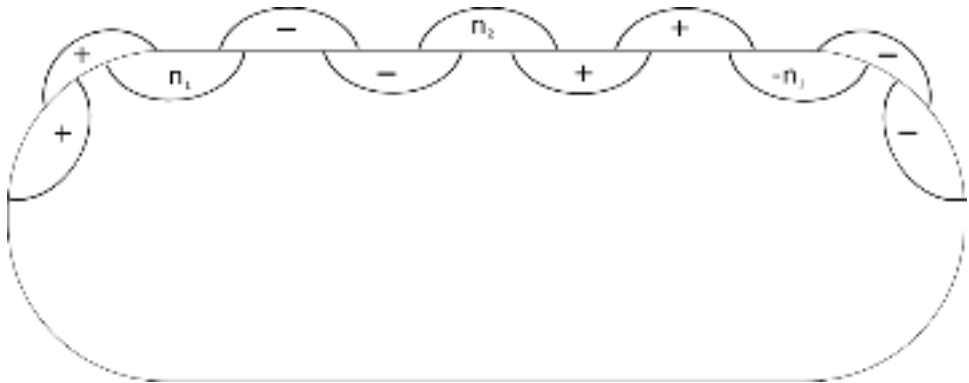

Fig. 15: Family of chord diagram of $\textbf{Sem}_1$.

A knot from the Milnor-Fox family with coefficient $n_1 = 2$ is a knot $8_{20}$ and it has a trivial second elementary ideal. Consider its duplicate. The chord diagram of this pairing is obtained from the figure if we take the following values of the coefficients: $|n_1| = 2$ and $n_2 = 2$. The Alexander matrix of the knot, which is a duplicate of $8_{20}$ looks like:

$$\begin{pmatrix} (t^2 - t + 1)^2 & t \\ 0 & (t^2 - t + 1)^2 \end{pmatrix}$$

Note that it will also have a trivial second elementary ideal. Note that you can duplicate any knot from the $\textbf{Sem}_1$ family in the same way as it was done earlier.

The result is a new family $\textbf{Sem}_3$. It's clear that a representative of this family with coefficients $|n_1| = 2$, $n_2 = 2$, $n_1 = 2$ will have the trivial second ideal of Alexander. This process can be continued further. Thus, there are infinitely many bipartite knots with a trivial second Alexander ideal. The fact that all these knots are different follows from the fact that they have a different first elementary ideal □

# Bipartite knot and its chord diagram

## Chord diagrams and corresponding representation matrix

Let us describe the method of constructing the matrix $Q_D$ of the matched diagram $D$ from its chord diagram $C_D$. Let us carry out a detailed reasoning for negative chords (for positive chords they are completely analogous). First, consider the case of an inner chord.

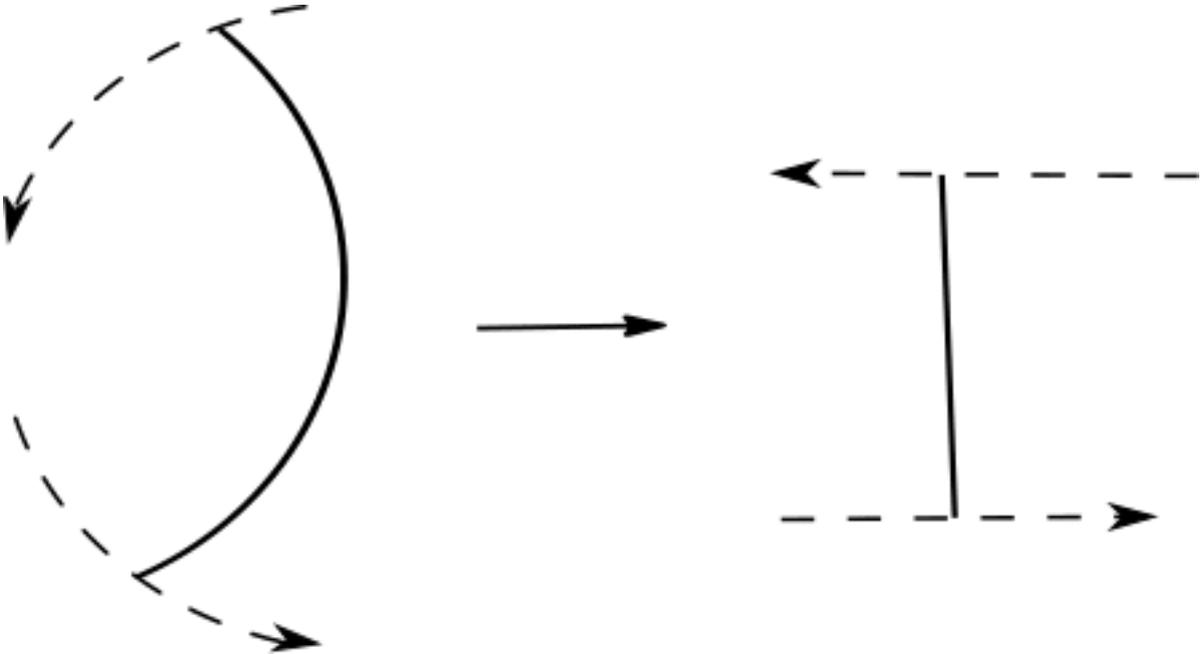

Fig. 16: Inner chord, local fragment of diagram

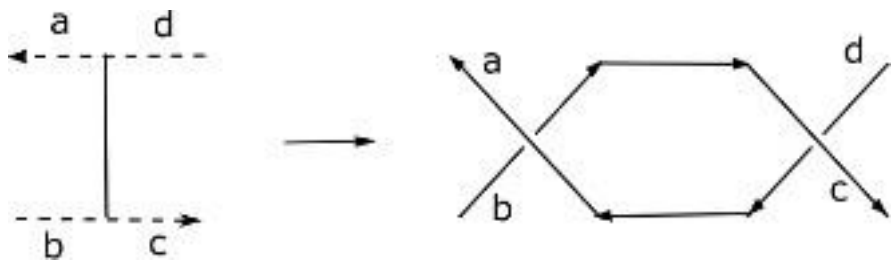

Fig. 17: Inner chord

Each chord of the chord diagram $C_D$ corresponds to a special pair of intersections.

Thus each chord corresponds to two relations in the fundamental group, and therefore a pair of rows in the matrix $Q_D$. The arcs of the diagram are in one-to-one correspondence with the arcs of the knot, and, therefore, with the generators of the fundamental group, which means that they correspond to the columns in the matrix. In this case $r_1 = aba^{-1}c^{-1}$ corresponds to the left crossing, and $r_2 = cdc^{-1}a^{-1}$ to the right (Fig. 17). Lets differentiate these relations according to the rules from the definition and then all group generators to the formal variable $t$.

1. $\frac{\partial r_1}{\partial a} \rightarrow 1 - t \qquad \frac{\partial r_2}{\partial a} \rightarrow -1$

2. $\frac{\partial r_1}{\partial b} \rightarrow t \qquad \frac{\partial r_2}{\partial b} \rightarrow 0$

3. $\frac{\partial r_1}{\partial c} \rightarrow -1 \qquad \frac{\partial r_2}{\partial c} \rightarrow 1 - t$

4. $\frac{\partial r_1}{\partial d} \rightarrow 0 \qquad \frac{\partial r_2}{\partial d} \rightarrow t$

It is convenient to describe the rule for filling the matrix in the following form:

you can schematically depict a chord with adjacent arcs. Near each arc will be written those elements of the matrix that are in the column corresponding to this arc, in the rows corresponding to this chord. Since one chord corresponds to two relations and two rows in the matrix, then each arc will correspond to two elements. Let's write inside the fragment of the diagram the elements corresponding to one relation and line, and outside to another. We get the following rule that allows you to write out the matrix Now consider the case of an outer chord, the special pair of crossings corresponding to the chord will have a different orientation (clockwise). (This is because the circle of the chord diagram is oriented counterclockwise, so the fragments of the diagram with outer chords must have a different direction to match the orientation of the diagram. (See Fig.19)). Relations for the case of an outer chord: $r_1 = aca^{-1}b^{-1}$ corresponds to the left crossing, and $r_2 = cac^{-1}d^{-1}$ to the right (Fig. 20). Now you can do the same steps as last time. The result is the following rule (Fig.21): Figures **??**

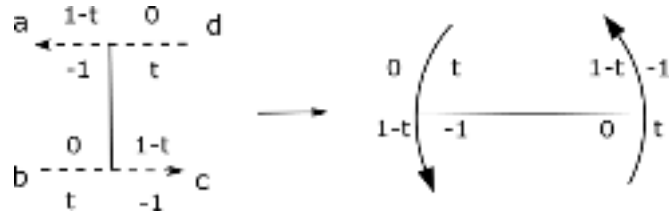

Fig. 18: Rule for the inner chord

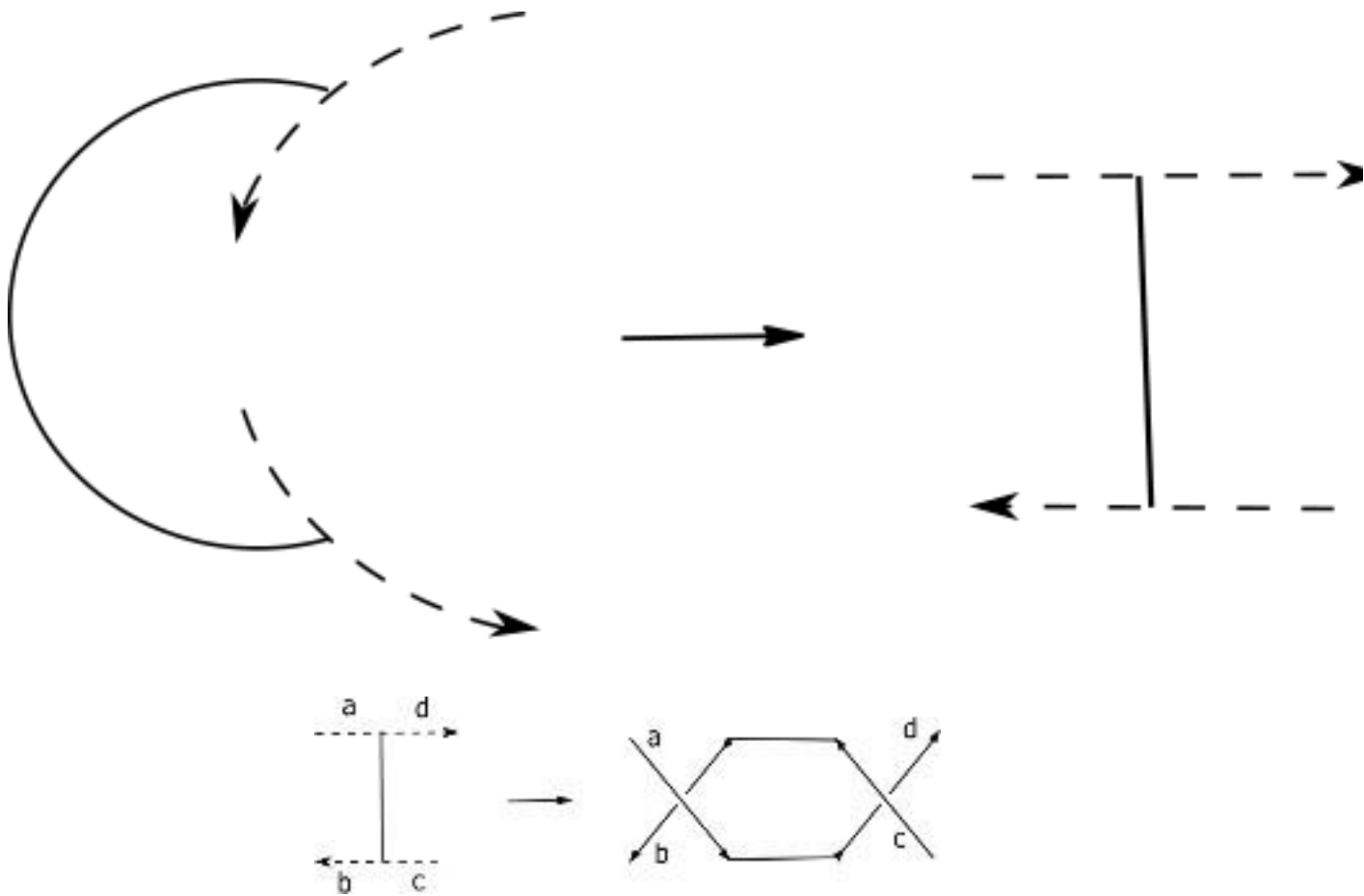

Fig. 20: Outer chord

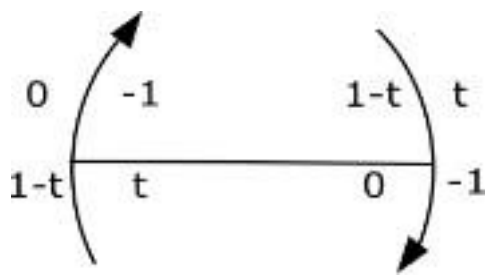

Fig. 21: .Rule for the outer chord

and **??** show all cases, on the left - outer chords, on the right - inner chords.

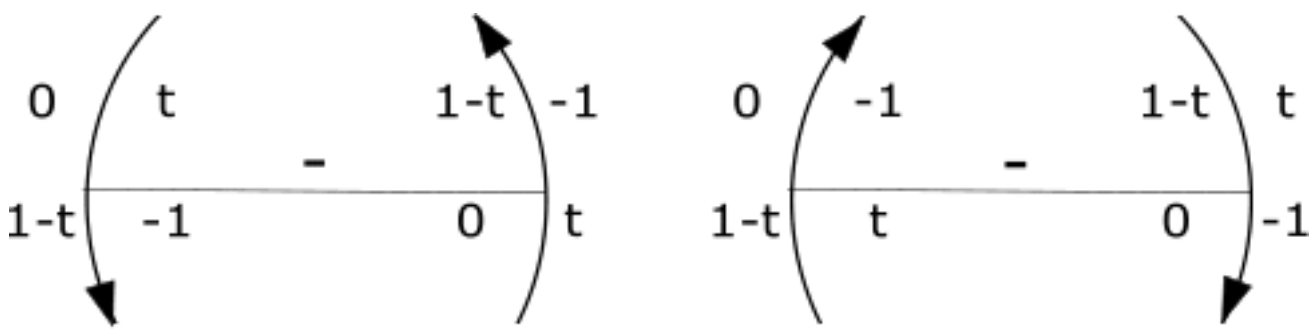

Fig. 22: Minuses

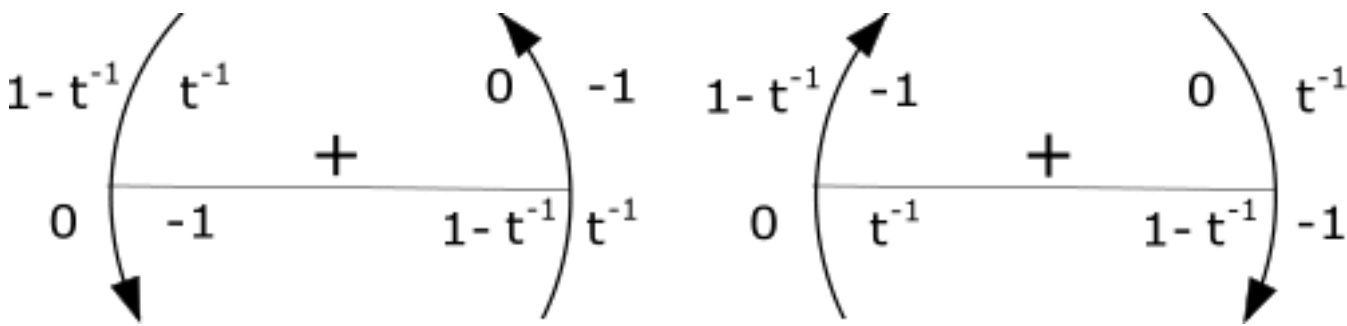

Fig. 23: Pluses

## Intersection graph and chord operations

It's clear that a diagram consisting exclusively of inner or outer chords corresponds to a trivial knot. We will call such chord diagrams trivial.

Any nontrivial chord diagram $C_D$ can be reduced to at least one of the following diagrams by removing a certain number of chords.

**Definition.** The intersection graph $G(D)$ of the chord diagram $C_D$ is a graph whose vertices correspond to the chords $C_D$, moreover, edges connect those and only those vertices for which the corresponding chords intersect if they are drawn as segments inside a circle.

The intersection graph is bipartite, and its vertices are marked in white and black in accordance with the belonging to a certain part (inner and outer chords). Also, each vertex of the graph has a sign depending on the sign of the chord to which it corresponds. Here is an example of building a graph using a chord diagram (Fig. 25).

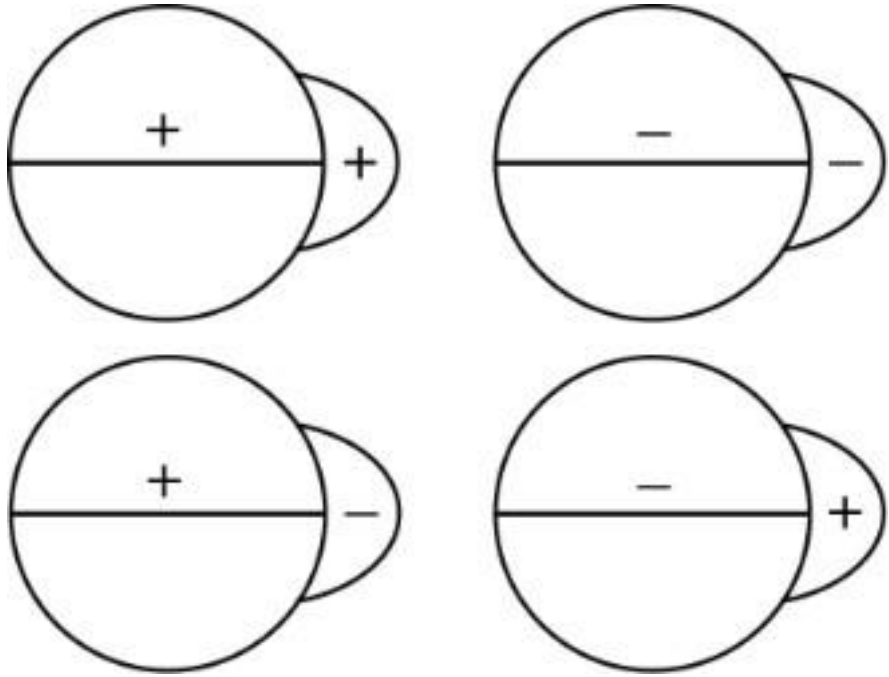

Fig. 24: Elementary diagrams

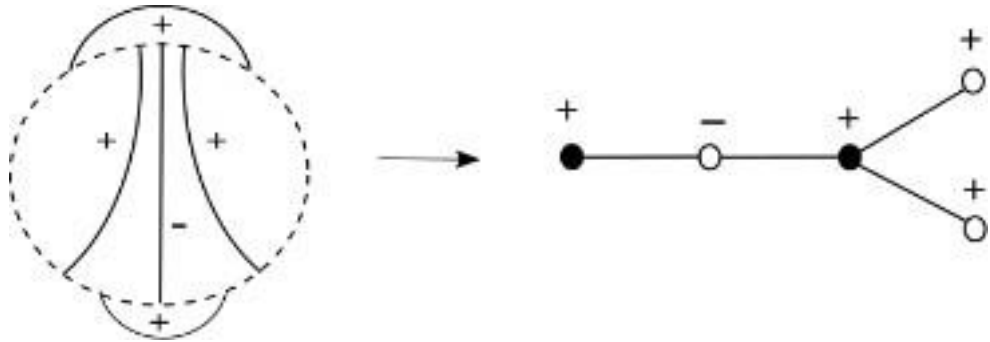

Fig. 25: Intersection graph for a chord diagram of a knot $8_{14}$.

We introduce two operations on chord diagrams, each of which increases the number of chords in the diagram (allowing us to obtain a new chord diagram).

**First operation:**

Let's choose some chord in the $C_D$ diagram. Let's add a new chord to the diagram so that both ends of the new chord are adjacent to one of the ends of the chord of the selected chord. (See Fig. 26 and 27) Obviously, the new chord will have two common arcs with the chosen one. Also, when adding a new chord, two new arcs appear on the diagram, we will assume that the new arcs are common arcs. In the figures **??** and **??**, the new chord is marked in green, and the new arcs are indicated by the letters $x$ and $y$.

In terms of the intersection graph, this means that, one new vertex is added to the intersection graph $G(D)$, which is connected by an edge to one of the vertices of the graph. The new vertex is indicated in green for illustration, although it is actually white if we add an inner chord and black if we add an outer one. Since the graph is connected, the vertex with which we connect a new one with an edge must be connected to at least one vertex from the graph, so at least one edge emanates from it, but it can be connected to several vertices, then it will already be several edges, these edges we marked with a dotted line.

**\* Second operation:**

Let's choose some chord in the $C_D$ diagram. Let's add a new chord to the diagram so that the ends of the new chord are adjacent to both ends of the selected chord. (See Fig. **??** - **??**) Obviously, the new chord will have two common arcs with the selected chord. Also, when adding a new chord, two new arcs appear on the diagram, we will assume that the new arcs are common arcs. In the figures **??** and **??**, the new chord is marked in green, and the new arcs are indicated by the letters $x$ and $y$.

In terms of the intersection graph, this means that a new vertex is added to the intersection graph, which is connected to all the vertices to which one of the graph vertices is connected. (See Fig. **??** - **??**) The new vertex is marked in green.

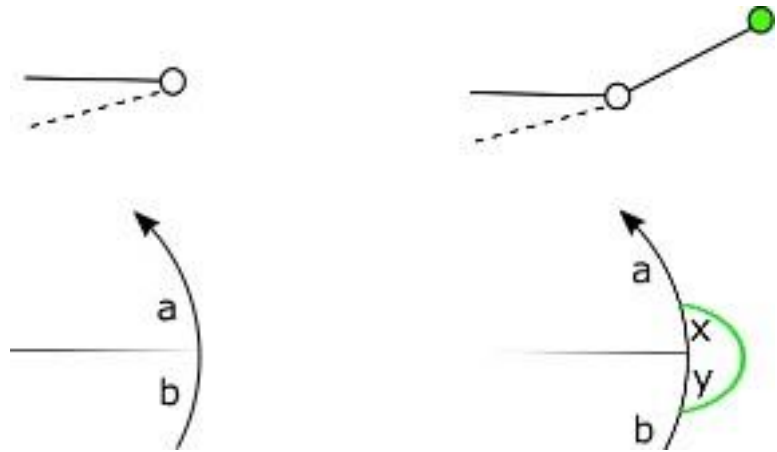

Fig. 26: First operation:adding an outer chord

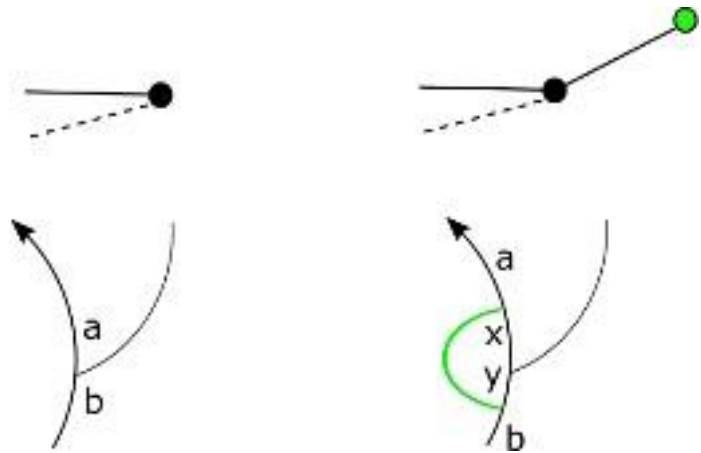

Fig. 27: First operation:adding an inner chord

*Remark.* Obviously, not all chord diagrams can be obtained by adding chords only in these ways. For example, a chord diagram corresponding to the knot $12n881$ (see Fig. 30) is impossible.

## Support chords and the second elementary ideal

**Definition.** A support chord is such a chord in a chord diagram $C_D$ with $k$ chords such that when deleting the rows of the matrix corresponding to it, among the minors of order $2k - 2$, there is an invertible one.

*Remark.* Obviously, if the chord diagram of a bipartite knot $K$ contains a support chord, then the second elementary ideal of this knot is trivial. In the simplest nontrivial chord diagrams with two chords, each chord is such.

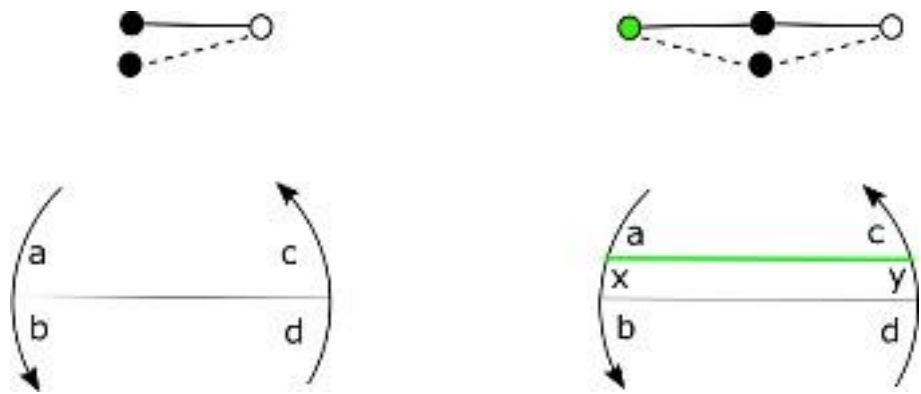

Fig. 28: Second operation:adding an inner chord

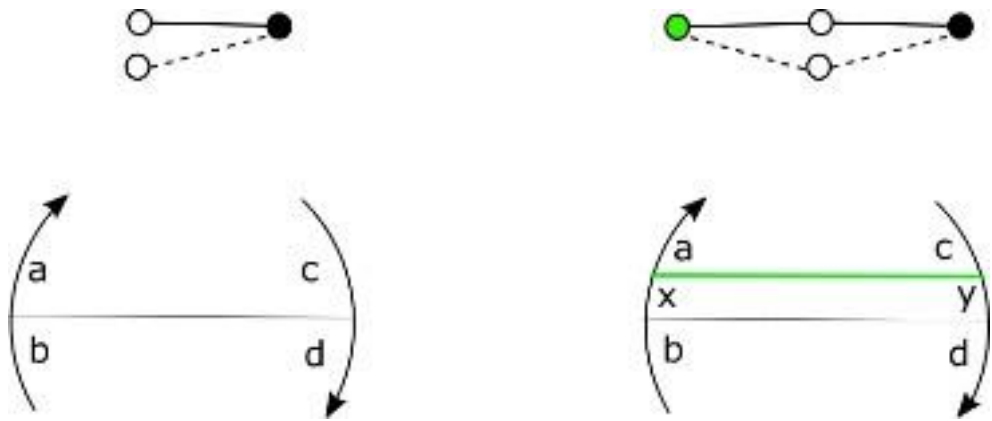

Fig. 29: Second operation:adding an outer chord

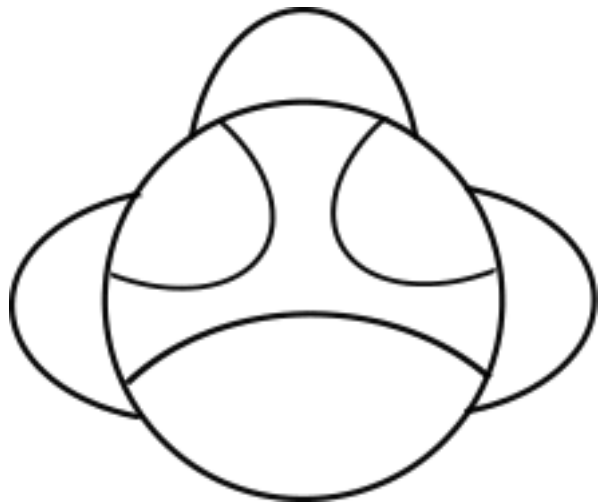

Fig. 30: $12n881$ without signs

**Lemma 5.** *Let $C_D$ be a diagram with a support chord $\chi_1$, then when a new chord is added to it using the two operations described above, $\chi_1$ stay support. Moreover, when adding a chord using the second operation, the new chord will also be support.*

Let the diagram $C_D$ correspond to the matrix $A_1$.

* **First operation:** We will assume that $\chi_1$ is positive. In the case when it is negative, all reasoning is similar.

Let's add a new chord $\chi_2$ using the first operation. There are four possible cases: inner positive chord, inner negative, outer positive, outer negative. We will limit to analyzing one of four cases. We will assume that we have added an inner negative chord. The reader, if desired, can apply the same reasoning in other cases and make sure that they are similar to the one discussed. Consider the $A_2$ matrix that corresponds to the new diagram. Without loss of generality, we will assume that the columns corresponding to the two new arcs of the diagram are the last, and the rows corresponding to the new chord $\chi_2$ are also. In the matrix $A_1$, we assume that the columns corresponding to the pair of arcs $\chi_1$ that become the arcs $\chi_2$ are the last, and the rows corresponding to $\chi_1$ are first.

**Lemma 6.** *The $A_2$ matrix is equivalent to the following matrix.*

$$
\begin{array}{c}
\chi_1 \\
\\
\\
\chi_2
\end{array}
\begin{array}{cc}
& \begin{array}{cc} x & y \end{array} \\
\left[
\begin{array}{cc|cc}
 & & 1-t & 0 \\
\quad A_1 & & t & -1 \\
 & & \vdots & \vdots \\
\hline
\vdots \ \vdots \ \ 1-t \ \ t-1 & & 1-t & t \\
\vdots \ \vdots \ \ t-1 \ \ 1-t & & -1 & 0
\end{array}
\right]
\end{array}
$$

Fig. 31

*Remark.* Dots are implied wherever zeros appears.

*Proof.* Let's add a pair of columns corresponding to new arcs to a pair of columns with nonzero elements in the last two rows. This obviously follows from the filling matrix rules illustrated above. □

**Example 3.** *Let's look at an example. Let's take one of the elementary diagrams $D_1$ (It is shown on the left in Fig. 32).*

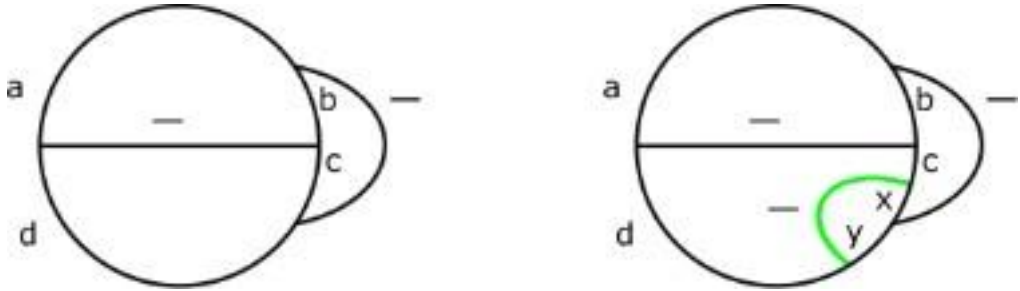

Fig. 32: $D_1$ and $D_2$.

*It corresponds to the matrix*

$$
\begin{array}{cccc}
a & b & c & d
\end{array}
$$
$$
\begin{bmatrix}
-1 & t & 0 & 1-t \\
0 & 1-t & -1 & t \\
t & 1-t & 0 & -1 \\
0 & -1 & t & 1-t
\end{bmatrix}
$$

Fig. 33

*The $D_2$ diagram is obtained from $D_1$ by adding a chord. It corresponds to the matrix If we add the columns corresponding to the arcs x y to the columns corresponding to the arcs c d, we get the matrix*

$$\begin{array}{c} \\ \chi_1 \\ \\ \\ \\ \chi_2 \\ \\ \end{array} \begin{array}{cccccc} a & b & c & d & x & y \\ \left[\begin{array}{cccccc} -1 & t & 0 & 0 & 0 & 1-t \\ 0 & 1-t & 0 & 0 & -1 & t \\ t & 1-t & 0 & -1 & 0 & 0 \\ 0 & -1 & t & 1-t & 0 & 0 \\ 0 & 0 & 1-t & t & 0 & -1 \\ 0 & 0 & -1 & 0 & t & 1-t \end{array}\right] \end{array}$$

Fig. 34

$$\begin{array}{c} \\ \chi_1 \\ \\ \\ \\ \chi_2 \\ \\ \end{array} \begin{array}{cccccc} a & b & c & d & x & y \\ \left[\begin{array}{cccccc} -1 & t & 0 & 1-t & 0 & 1-t \\ 0 & 1-t & -1 & t & -1 & t \\ t & 1-t & 0 & -1 & 0 & 0 \\ 0 & -1 & t & 1-t & 0 & 0 \\ 0 & 0 & 1-t & t-1 & 0 & -1 \\ 0 & 0 & t-1 & 1-t & t & 1-t \end{array}\right] \end{array}$$

Fig. 35

It obviously follows from the lemma and the conditions of the theorem that among the minors of order $n-2$ obtained by deleting the rows corresponding to the chord $\chi_1$, there is trivial. This is clear, since when deleting the required lines in the last two In columns, nonzero elements remain only in the last two rows, while the determinant of the 22 matrix in the lower right corner is equal to the unit of the ring.

Second operation: Let's assume that $\chi$ is negative inner chord. In the rest cases, all reasoning is similar. Now add a new chord $\chi_2$ using second operation. In order to prove the statement, we need the following. lemma.

**Lemma 7.** *The $A_2$ matrix is equivalent to the following matrix*

Fig. 36

*Proof.* Let's add a pair of columns corresponding to new arcs to a pair of columns with nonzero elements in the last two rows. This obviously follows from the rules for filling the matrix illustrated above. □

**Example 4.** *Let's look at an example. Let's take one of the elementary diagrams $D_1$*
*(It is shown on the left in Fig. 37).*

It corresponds to the matrix

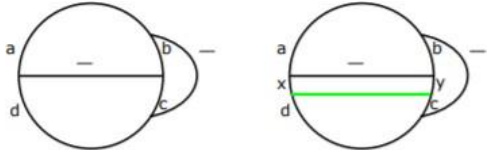

Fig. 37

$$
\begin{array}{cccc}
a & b & c & d \\
\end{array}
$$
$$
\begin{bmatrix}
t & 1\text{-}t & 0 & \text{-}1 \\
0 & \text{-}1 & t & 1\text{-}t \\
\text{-}1 & t & 0 & 1\text{-}t \\
0 & 1\text{-}t & \text{-}1 & t
\end{bmatrix}
$$

Fig. 38

The $D_2$ diagram is obtained from $D_1$ by adding a chord. It corresponds to the matrix

If we add the columns corresponding to the arcs $x\,y$ to the columns corresponding to the arcs $c\,d$, we get the matrix shown in Fig. 40. Based on the same considerations as in the previous paragraph, it is obvious that in this case the chord $chi_1$ stay support. But if we delete the lines corresponding to the chord $\chi_2$, then in the last two columns the only nonzero elements will remain only in the rows corresponding to the chord $\chi_1$. The determinant of the $2\times 2$ matrix in the upper right corner is equal to the unit of the ring. Hence, no element from the rows corresponding to $\chi_1$ in the submatrix $A_1$ will be included in the minors, thus the rows corresponding to $\chi_1$ were "deleted" from the matrix $A_1$. By the condition of the theorem, there is a trivial minor, and hence the chord $\chi_2$ is also a support chord.

## Rational knots

Let $p$ and $q$ be coprime integers, $p > 0$, $|\frac{p}{q}| < 1$ Consider the continued fraction expansion of the number $\frac{p}{q}$.

$$
\frac{p}{q} = \cfrac{1}{b_1 + \cfrac{1}{b_2 + \cfrac{1}{\ldots + \cfrac{1}{b_{n-1} + \cfrac{1}{b_n}}}}}
$$

where $b_i$ are nonzero integers. Further, we will use the shorter notation $[b_1, b_2, \ldots, b_{n-1}, b_n]$ to denote a continued fraction with denominators $b_1, b_2, \ldots, b_{n-1}, b_n$.

$$
\begin{array}{ccccccc}
 & a & b & c & d & x & y \\
\end{array}
$$
$$
\begin{array}{c}
\chi_1 \\ \\ \\ \\ \chi_2 \\ \\
\end{array}
\begin{bmatrix}
t & 1\text{-}t & 0 & 0 & \text{-}1 & 0 \\
0 & \text{-}1 & 0 & 0 & 1\text{-}t & t \\
\text{-}1 & t & 0 & 1\text{-}t & 0 & 0 \\
0 & 1\text{-}t & \text{-}1 & t & 0 & 0 \\
0 & 0 & t & 1\text{-}t & 0 & \text{-}1 \\
0 & 0 & 0 & \text{-}1 & t & 1\text{-}t
\end{bmatrix}
$$

Fig. 39

|  | $a$ | $b$ | $c$ | $d$ | $x$ | $y$ |
|---|---|---|---|---|---|---|
| $\chi_1$ | $t$ | $1-t$ | $0$ | $-1$ | $-1$ | $0$ |
|  | $0$ | $-1$ | $t$ | $1-t$ | $1-t$ | $t$ |
|  | $-1$ | $t$ | $0$ | $1-t$ | $0$ | $0$ |
|  | $0$ | $1-t$ | $-1$ | $t$ | $0$ | $0$ |
| $\chi_2$ | $0$ | $0$ | $t-1$ | $1-t$ | $0$ | $-1$ |
|  | $0$ | $0$ | $1-t$ | $t-1$ | $t$ | $1-t$ |

Fig. 40

Consider a braid on four strands corresponding to the word $A^{b_1}B^{b_2}A^{b_3}\ldots$, where $A$ and $B$ are the fragments shown below and connect from left to right. Then, take the closure of this braid depending on the parity of $n$.

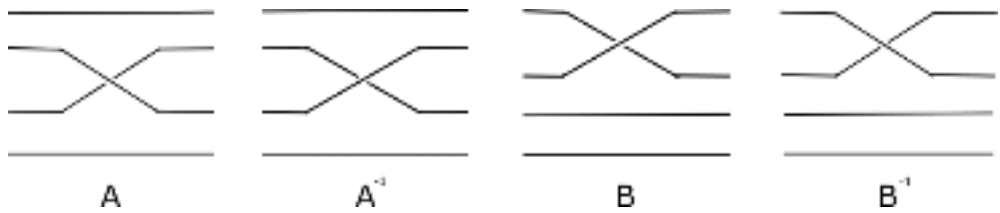

Fig. 41: Fragments

The knot (link) diagram constructed in this way will be called the natural diagram of the rational knot (link)

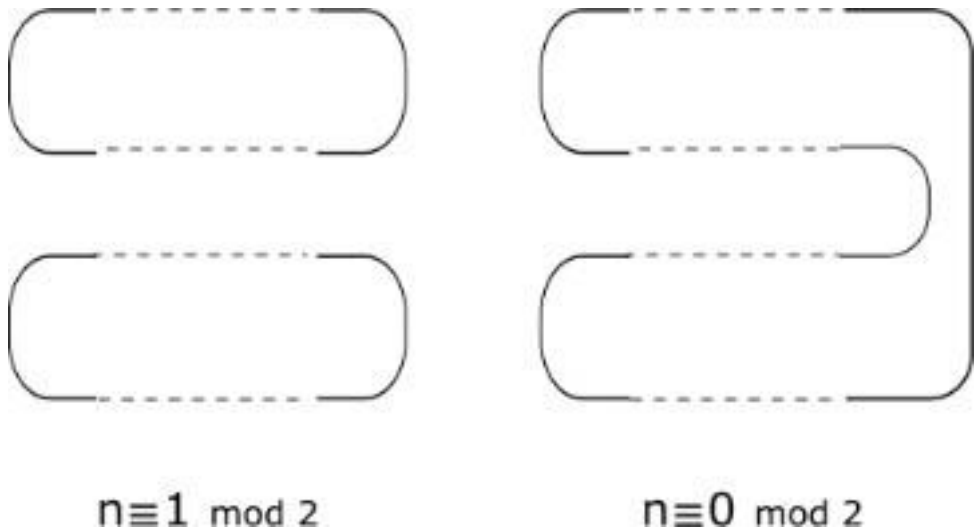

Fig. 42: Closure

$L(\frac{p}{q})$. Any rational knot has a matched diagram, since any rational number can be represented as a continued fraction with even denominators (positive or negative) Thus, the continued fraction $[2a_1, 2a_2, \ldots, 2a_{2n-1}, 2a_{2n}]$ corresponds to a chord diagram with sequentially added $a_{2n}$ inner chords, $a_{2n-1}$ external, $\ldots$, $a_1$ by outer chords.

In what follows we will write $[a_1, \ldots, a_{2n}]$ instead of $[2a_1, 2a_2, \ldots, 2a_{2n-1}, 2a_{2n}]$

**Lemma 8.** *The chord diagram of any rational knot has a support chord.*

*Proof.* We carry out the proof by induction on the number of denominators in the continued fraction expansion of a rational number. Induction base. Note that for any rational knot diagram corresponding to the sequence $[a_1, a_2]$ with two members, the last added chord is supporting. This is true, since you can choose an elementary diagram, the sign of the inner chord of which corresponds to the sign $a_2$, and the outer one corresponds to the sign $a_1$, then add $a_2 - 1$ to the outer chord of the elementary diagram inner chords using the first operation. By the lemma proved earlier, the outer chord stay support. And then add $a_1 - 1$ of parallel outer chords to the outer chord. By the lemma proved earlier, all outer chords will be supporting, in particular, the last added one. (See Fig. 45) The figure shows the steps described step by step. At each step, the added chords are indicated in green.

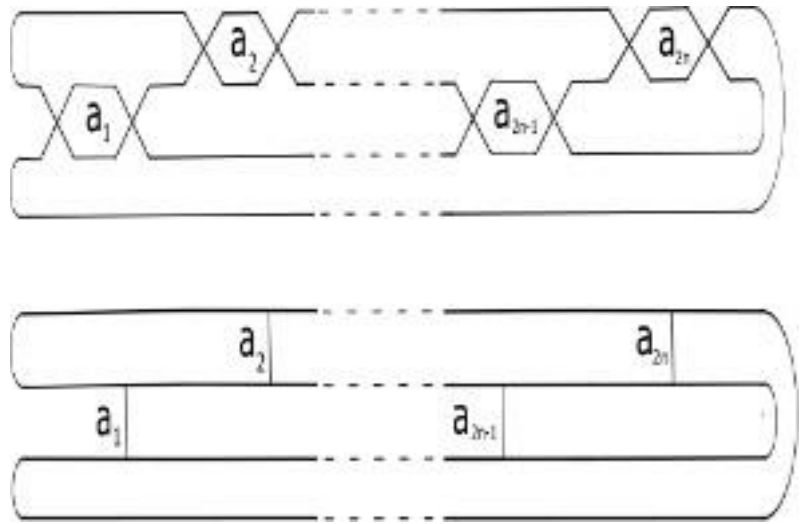

Fig. 43: Rational knot and its chord diagram

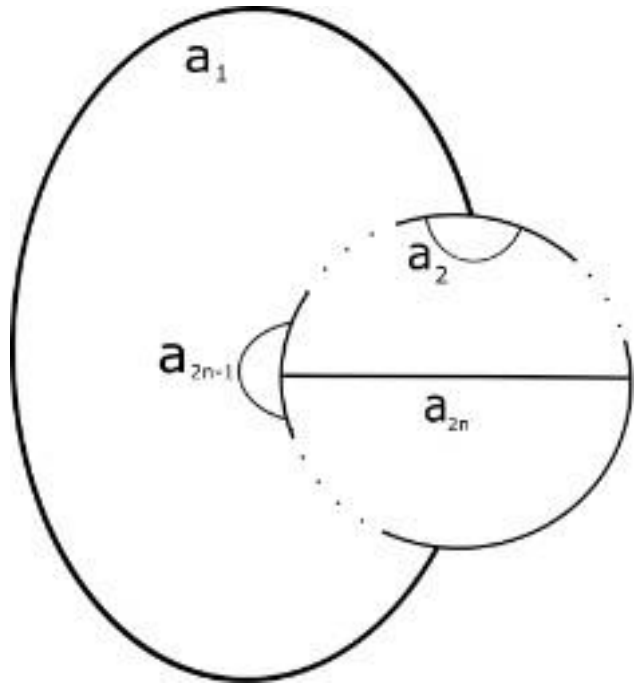

Fig. 44: Chord diagram of rational knot

Let this be true for chord diagrams of rational knots corresponding to the expansion with $2n$ denominators.

Any sequence of even numbers with $2n + 2$ members $[b_1, ...b_{2n+2}]$ corresponds to the sequence $[a_1, ..., a_{2n}]$ such such that $a_i = b_{i+2}$, $1 \le i \le 2n$

Take the chord diagram of the knot corresponding to the sequence $[b_1, ...b_{2n+2}]$. Let us show that it can be constructed from the diagram corresponding to $[a_1, ...., a_{2n}]$. Consider the diagram corresponding to $[a_1, \quad , a_{2n}]$. Add to the last added outer chord in parallel one more outer chord $\chi$ of the sign, as in the number $b_1$. The new chord $\chi$ is support. Let us add to $\chi$ $b_2$ internal sign chords of the number $b_2$ by the first operation, $\chi$ stay support. Add to $\chi$ $b_1 \_ 1$ parallel sign chords of the number $b_1$. Each of these added outer chords is support, including the last one added. Note that we have obtained a chord diagram of the knot corresponding to $[b_1, ...b_{2n+2}]$. (See Fig. 46) The figure shows the steps described step by step. At each step, the added chords are indicated in green.

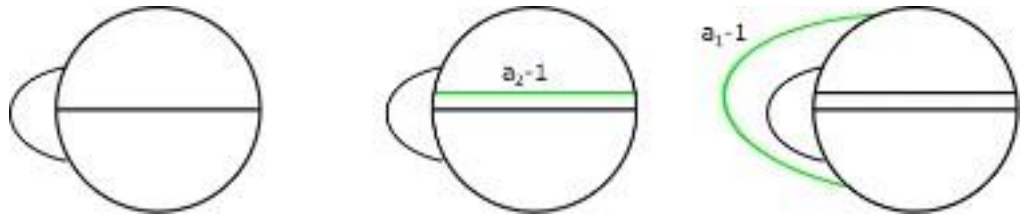

Fig. 45: Chord diagrams

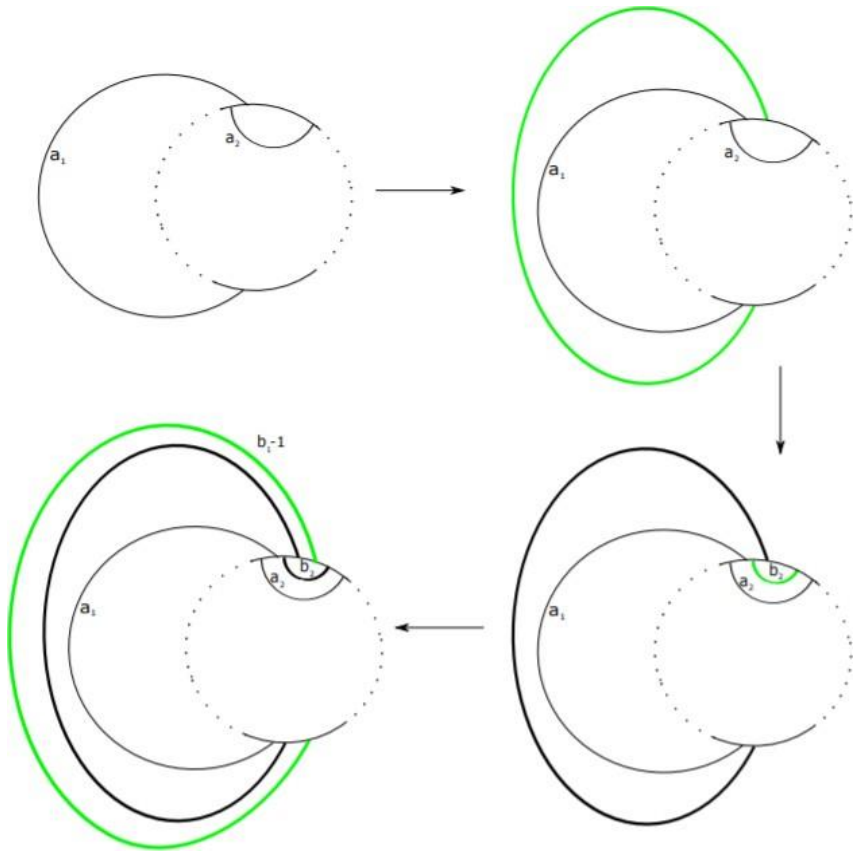

Fig. 46: Chord diagrams